\documentclass[10pt,twoside]{article}
\usepackage{Latex-document,amssymb}
\usepackage{multicol}
\usepackage{graphicx}
\def\volumeno{I}
\title{\bf The Work of Laurent Lafforgue\vskip 6mm}
\author{G\'{e}rard Laumon\thanks{CNRS and
Universit\'{e} Paris-Sud, UMR 8628, Math\'{e}matique, F-91405 Orsay
Cedex, France, gerard.laumon@math.u-psud.fr}\vspace*{-0.5cm}}
\date{\vspace{-8mm}}

\begin{document}

\maketitle \thispagestyle{first} \setcounter{page}{91}

\vskip 12mm

Laurent Lafforgue has been awarded the Fields Medal for his proof of the
Langlands correspondence for the full linear groups $\mathop{\rm
GL}\nolimits_{r}$ ($r\geq 1$) over function fields.

What follows is a brief introduction to the Langlands correspondence and
to Lafforgue's theorem.

\section{The Langlands correspondence}

\vskip-5mm \hspace{5mm }

A global field is either a number field, i.e. a finite extension of
$\mathbb{Q}$, or a function field of characteristic $p>0$ for some prime
number $p$, i.e. a finite extension of $\mathbb{F}_{p}(t)$ where
$\mathbb{F}_{p}$ is the finite field with $p$ elements.  The global
fields constitute a primary object of study in number theory and
arithmetic algebraic geometry.

The conjectural Langlands correspondence, which was first formulated by
Robert Langlands in 1967 in a letter to Andr\'{e} Weil, relates two
fundamental objects which are naturally attached to a global field $F$:
\vskip 1mm

\begin{description}
\item{\hspace*{3mm} --} its Galois group $\mathop{\rm Gal}(\overline{F}/F)$,
where $\overline{F}$ is an algebraic closure of $F$, or more accurately
its motivic Galois group of $F$ which is by definition the tannakian
group of the tensor category of Grothendieck motives over $F$, \vskip
1mm

\item{\hspace*{3mm} --} the ring $\mathbb{A}$ of ad\`{e}les of $F$, or more
precisely the collection of Hilbert spaces $L^{2}(G(F)\backslash
G(\mathbb{A}))$ for all reductive groups $G$ over $F$.
\end{description}

Roughly speaking, for any (connected) reductive group $G$ over $F$,
Langlands introduced a dual group ${}^{{\rm L}}G=\widehat{G}\rtimes
\mathop{\rm Gal}(\overline{F}/F)$, the connected component $\widehat{G}$
of which is the complex reductive group whose roots are the co-roots of
$G$ and vice versa.  And he predicted that a large part of the spectral
decomposition of the Hilbert space $L^{2}(G(F)\backslash
G(\mathbb{A}))$, equipped with the action by right translations of
$G(\mathbb{A})$, is governed by representations of the motivic Galois
group of $F$ with values in ${}^{{\rm L}}G$.

Of special importance is the group $G=\mathop{\rm GL}_{r}$, the
Langlands dual of which is simply the direct product ${}^{{\rm L}}
\mathop{\rm GL}\nolimits_{r}=\mathop{\rm GL}\nolimits_{r}(\mathbb{C})
\times \mathop{\rm Gal} (\overline{F}/F)$.  Indeed, any complex
reductive group $\widehat{G}$ may be embedded into $\mathop{\rm
GL}\nolimits_{r}(\mathbb{C})$ for some $r$.

The particular case $G=\mathop{\rm GL}_{1}$ of the Langlands
correspondence is the abelian class field theory of Teiji Takagi and
Emil Artin which was developed in the 1920s as a wide extension of the
quadratic reciprocity law.

The Langlands correspondence embodies a large part of number theory,
arithmetic algebraic geometry and representation theory of Lie groups.
Small progress made towards this conjectural correspondence had already
amazing consequences, the most striking of them being the proof of
Fermat's last theorem by Andrew Wiles.  Famous conjectures, such as the
Artin conjecture on $L$-functions and the Ramanujan-Petersson
conjecture, would follow from the Langlands correspondence.

\section{Lafforgue's main theorem}

\vskip-5mm \hspace{5mm }

Over number fields, the Langlands correspondence in its full generality
seems still to be out of reach.  Even its precise formulation is very
involved.  In the function field case the situation is much better.
Thanks to Lafforgue, the Langlands correspondence for $G=\mathop{\rm
GL}\nolimits_{r}$ is now completely understood.

 From now on, $F$ is a function field of characteristic $p>0$. We also
fix some auxiliary prime number $\ell\not=p$.

As Alexandre Grothendieck showed, any algebraic variety over $F$ gives
rise to $\ell$-adic representations of $\mathop{\rm
Gal}(\overline{F}/F)$ on its \'{e}tale cohomology groups and the
irreducible $\ell$-adic representations of $\mathop{\rm
Gal}(\overline{F}/F)$ are good substitutes for irreducible motives over
$F$.  Therefore, the Langlands correspondence may be nicely formulated
using $\ell$-adic representations.

Let $r$ be a positive integer.  On the one hand, we have the set ${\cal
G}_{r}$ of isomorphism classes of rank $r$ irreducible $\ell$-adic
representations of $\mathop{\rm Gal} (\overline{F}/F)$ the determinant
of which is of finite order.  To each $\sigma\in {\cal G}_{r}$,
Grothendieck attached an Eulerian product $L(\sigma
,s)=\prod_{x}L_{x}(\sigma ,s)$ over all the places $x$ of $F$, which is
in fact a rational function of $p^{-s}$ and which satisfies a functional
equation of the form $L(\sigma ,s)=\varepsilon (\sigma
,s)L(\sigma^{\vee},1-s)$ where $\sigma^{\vee}$ is the contragredient
representation of $\sigma$ and $\varepsilon (\sigma ,s)$ is some
monomial in $p^{-s}$.  If $\sigma$ is unramified at a place $x$, we have
$$
L_{x}(\sigma ,s)=\prod_{i=1}^{r}{1\over 1-z_{i}p^{-s\mathop{\rm
deg}(x)}}
$$
where $z_{1},\ldots ,z_{r}$ are the {\it Frobenius eigenvalues} of
$\sigma$ at $x$ and $\mathop{\rm deg}(x)$ is the degree of the place
$x$.

On the other hand, we have the set ${\cal A}_{r}$ of isomorphism classes
of cuspidal automorphic representations of $\mathop{\rm
GL}_{r}(\mathbb{A})$ the central character of which is of finite order.
Thanks to Langlands' theory of Eisenstein series, they are the building
blocks of the spectral decomposition of $L^{2}(\mathop{\rm
GL}\nolimits_{r}(F)\backslash \mathop{\rm
GL}\nolimits_{r}(\mathbb{A}))$. To each $\pi\in {\cal A}_{r}$, Roger
Godement and Herv\'{e} Jacquet attached an Eulerian product $L(\pi
,s)=\prod_{x} L_{x}(\pi ,s)$ over all the places $x$ of $F$, which is
again a rational function of $p^{-s}$, satisfying a functional equation
$L(\pi ,s)=\varepsilon (\pi ,s)L(\pi^{\vee},1-s)$.  If $\pi$ is
unramified at a place $x$, we have
$$
L_{x}(\pi ,s)=\prod_{i=1}^{r}{1\over 1-z_{i}p^{-s\mathop{\rm deg}(x)}}
$$
where $z_{1},\ldots ,z_{r}$ are called the {\it Hecke eigenvalues} of
$\pi$ at $x$. \vskip 1mm

{\bf Theorem} \it {\rm (i)} {\rm (The Langlands Conjecture\rm )} There
is a unique bijective correspondence $\pi\rightarrow \sigma (\pi )$,
preserving $L$-functions in the sense that $L_{x}(\sigma (\pi ),s)=
L_{x}(\pi ,s)$ for every place $x$, between ${\cal A}_{r}$ and ${\cal
G}_{r}$.

{\rm (ii)} {\rm (The Ramanujan-Petersson Conjecture\rm )} For any
$\pi\in {\cal A}_{r}$ and for any place $x$ of $F$ where $\pi$ is
unramified, the Hecke eigenvalues $z_{1},\ldots ,z_{r}\in
\mathbb{C}^{\times}$ of $\pi$ at $x$ are all of absolute value $1$.

{\rm (iii)} {\rm (The Deligne Conjecture\rm )} Any $\sigma\in {\cal
G}_{r}$ is pure of weight zero, i.e. for any place $x$ of $F$ where
$\sigma$ is unramified, and for any field embedding $\iota :
\overline{\mathbb{Q}}_{\ell}\hookrightarrow \mathbb{C}$, the images
$\iota (z_{1}),\ldots ,\iota (z_{r})$ of the Frobenius eigenvalues of
$\sigma$ at $x$ are all of absolute value $1$. \rm \vskip 1mm

As I said earlier, in rank $r=1$, the theorem is a reformulation of the
abelian class field theory in the function field case.  Indeed, the
reciprocity law may be viewed as an injective homomorphism with dense
image
$$
F^{\times}\backslash \mathbb{A}^{\times}\rightarrow \mathop{\rm
Gal}(\overline{F}/F)^{{\rm ab}}
$$
from the id\`{e}le class group to the maximal abelian quotient of the
Galois group.

In higher ranks $r$, the first breakthrough was made by Vladimir
Drinfeld in the 1970s.  Introducing the fundamental concept of shtuka,
he proved the rank $r=2$ case.  It is a masterpiece for which, among
others works, he was awarded the Fields Medal in 1990.

\section{The strategy}

\vskip-5mm \hspace{5mm }

The strategy that Lafforgue is following, and most of the geometric
objets that he is using, are due to Drinfeld.  However, the gap between
the rank two case and the general case was so big that it took more than
twenty years to fill it.

Lafforgue considers the $\ell$-adic cohomology of the moduli stack of
rank $r$ Drinfeld shtukas (see the next section) as a representation of
$\mathop{\rm GL}\nolimits_{r}(\mathbb{A}) \times \mathop{\rm
Gal}(\overline{F}/F)\times \mathop{\rm Gal} (\overline{F}/F)$.  By
comparing the Grothendieck-Lefschetz trace formula (for Hecke operators
twisted by powers of Frobenius endomorphisms) with the Arthur-Selberg
trace formula, he tries to isolate inside this representation a
subquotient which decomposes as
$$
\bigoplus_{\pi\in {\cal A}_{r}}\pi\otimes\sigma (\pi )^{\vee}\otimes
\sigma (\pi ).
$$

Such a comparison of trace formulas was first made by Yasutaka Ihara in
1967 for modular curves over $\mathbb{Q}$.  Since, it has been
extensively used for Shimura varieties and Drinfeld modular varieties by
Langlands, Robert Kottwitz and many others.  There are two main
difficulties to overcome to complete the comparison: \vskip 1mm

\begin{description}
\item{\hspace*{3mm} --} to prove suitable cases of a combinatorial conjecture of Langlands
and Diana Shelstad, which is known as the {\it Fundamental Lemma},
\vskip 1mm

\item{\hspace*{3mm} --} to compare the contribution of the ``fixed points at infinity'' in
the Grothendieck-Lefschetz trace formula with the weighted orbital
integrals of James Arthur which occur in the geometric side of the
Arthur-Selberg trace formula.
\end{description}

For the moduli space of shtukas, the required cases of the Fundamental
Lemma were proved by Drinfeld in the 1970s.  So, only the second
difficulty was remaining after Drinfeld had completed his proof of the
rank $2$ case.  This is precisely the problem that Lafforgue has solved
after seven years of very hard work.  The proof has been published in
three papers totalling about 600 pages.

\section{Drinfeld shtukas}

\vskip-5mm \hspace{5mm }

Let $X$ be ``the'' smooth, projective and connected curve over
$\mathbb{F}_{p}$ whose field of rational functions is $F$.  It plays the
role of the ring of integers of a number field.  Its closed points are
the places of $F$. For any such point $x$ we have the completion $F_{x}$
of $F$ at $x$ and its ring of integers ${\cal O}_{x}\subset F_{x}$.

Let ${\cal O}=\prod_{x}{\cal O}_{x}\subset \mathbb{A}$ be the maximal
compact subring of the ring of ad\`{e}les. Weil showed that the double
coset space
$$
\mathop{\rm GL}\nolimits_{r}(F)\backslash \mathop{\rm GL}
\nolimits_{r}(\mathbb{A})/\mathop{\rm GL}\nolimits_{r}({\cal O})
$$
can be naturally identified with the set of isomorphism classes of rank
$r$ vector bundles on $X$.

Starting from this observation, with the goal of realizing a congruence
relation between Hecke operators and Frobenius endomorphisms, Drinfeld
defined a rank $r$ {\it shtuka} over an arbitrary field $k$ of
characteristic $p$ as a diagram
$$
{}^{\tau}{\cal E}\,\smash{\mathop{\hbox to 8mm{\rightarrowfill}}
\limits_{\scriptstyle \varphi}^{\scriptstyle\sim}}\,\overbrace{{\cal
E}''\hookrightarrow {\cal E}'\hookleftarrow {\cal E}}^{{\rm Hecke}}
$$
where ${\cal E}$, ${\cal E}'$ and ${\cal E}'$ are rank $r$ vector
bundles on the curve $X_{k}$ deduced from $X$ by extending the scalars
to $k$, where ${\cal E}\hookrightarrow {\cal E}'$ is an elementary upper
modification of ${\cal E}$ at some $k$-rational point of $X$ which is
called the {\it pole} of the shtuka, where ${\cal E}'' \hookrightarrow
{\cal E}'$ is an elementary lower modification of ${\cal E}'$ at some
$k$-rational point of $X$ which is called the {\it zero} of the shtuka,
and where ${}^{\tau}{\cal E}$ is the pull-back of ${\cal E}$ by the
endomorphism of $X_{k}$ which is the identity on $X$ and the Frobenius
endomorphism on $k$.

Drinfeld proved that the above shtukas are the $k$-rational points of an
algebraic stack over $\mathbb{F}_{p}$ which is equipped with a
projection onto $X\times X$ given by the pole and the zero.  More
generally, he introduced level structures on rank $r$ shtukas and he
constructed an algebraic stack $\mathop{\rm Sht}\nolimits_{r}$
parametrizing rank $r$ shtukas equipped with a compatible system of
level structures.  This last algebraic stack is endowed with an
algebraic action of $\mathop{\rm GL}\nolimits_{r}(\mathbb{A})$ through
the Hecke operators. \eject

\section{Iterated shtukas}

\vskip-5mm \hspace{5mm }

The geometry at infinity of the moduli stack $\mathop{\rm Sht}
\nolimits_{r}$ is amazingly complicated.  The algebraic stack
$\mathop{\rm Sht} \nolimits_{r}$ is not of finite type and one needs to
{\it truncate} it to obtain manageable geometric objects.  Bounding the
Harder-Narasimhan polygon of a shtuka, Lafforgue defines a family of
open substacks $(\mathop{\rm Sht} \nolimits_{r}^{\leq P})_{P}$ which are
all of finite type and whose union is the whole moduli stack.  But in
doing so, he loses the action of the Hecke operators which do not
stabilize those open substacks.

In order to recover the action of the Hecke operators, Lafforgue
enlarges $\mathop{\rm Sht} \nolimits_{r}$ by allowing specific
degenerations of shtukas that he has called {\it iterated shtukas}.

More precisely, Lafforgue lets the isomorphism $\varphi:{}^{\tau}{\cal
E}\,\smash{\mathop{\hbox to 8mm{\rightarrowfill}}\limits^{\scriptstyle
\sim}}\,{\cal E}''$ appearing in the definition of a shtuka, degenerate
to a {\it complete homomorphism} ${}^{\tau}{\cal E}\Rightarrow {\cal
E}''$, i.e. a continuous family of complete homomorphisms between the
stalks of the vector bundles ${}^{\tau}{\cal E}$ and ${\cal E}''$.

Let me recall that a complete homomorphism $V\Rightarrow W$ between two
vector spaces of the same dimension $r$ is a point of the partial
compactification $\mathop{\rm H\widetilde{om}}(V,W)$ of $\mathop{\rm
Isom}(V,W)$ which is obtained by successively blowing up the
quasi-affine variety $\mathop{\rm Hom}(V,W)-\{0\}$ along its closed
subsets
$$
\{f\in \mathop{\rm Hom}(V,W)-\{0\}\mid \mathop{\rm rank}(f)\leq i\}
$$
for $i=1,\ldots ,r-1$.  If $V=W$ is the standard vector space of
dimension $r$, the quotient of $\mathop{\rm H\widetilde{om}}(V,W)$ by
the action of the homotheties is the Procesi-De Concini compactification
of $\mathop{\rm PGL}\nolimits_{r}$.

In particular, Lafforgue obtains a smooth compactification, with a
normal crossing divisor at infinity, of any truncated moduli stack of
shtukas without level structure.

\section{One key of the proof}

\vskip-5mm \hspace{5mm }

Lafforgue proves his main theorem by an elaborate induction on $r$.
Compared to Drinfeld's proof of the rank $2$ case, a very simple but
crucial novelty in Lafforgue's proof is the distinction in the
$\ell$-adic cohomology of $\mathop{\rm Sht}\nolimits_{r}$ between the
$r$-negligible part (the part where all the irreducible constituents as
Galois modules are of dimension $<r$) and the $r$-essential part (the
rest).  Lafforgue shows that the difference between the cohomology of
$\mathop{\rm Sht}\nolimits_{r}$ and the cohomology of any truncated
stack $\mathop{\rm Sht}\nolimits_{r}^{\leq P}$ is $r$-negligible.  He
also shows that the cohomology of the boundary of $\mathop{\rm
Sht}\nolimits_{r}^{\leq P}$ is $r$-negligible. Therefore, the
$r$-essential part, which is defined purely by considering the Galois
action and which is naturally endowed with an action of the Hecke
operators, occurs in the $\ell$-adic cohomology of any truncated moduli
stack $\mathop{\rm Sht}\nolimits_{r}^{\leq P}$ and also in their
compactifications.

At this point, Lafforgue makes an extensive use of the proofs by Richard
Pink and Kazuhiro Fujiwara of a conjecture of Deligne on the
Grothendieck-Lefschetz trace formula.

\section{Compactification of thin Schubert cells}

\vskip-5mm \hspace{5mm }

In proving the Langlands conjecture for functions fields, Lafforgue
tried to construct nice compactifications of the truncated moduli stacks
of shtukas with arbitrary level structures.  A natural way to do that is
to start with some nice compactifications of the quotients of
$\mathop{\rm PGL} \nolimits_{r}^{n+1}/\mathop{\rm PGL}\nolimits_{r}$ for
all integers $n\geq 1$, and to apply a procedure similar to the one
which leads to iterated shtukas.

Lafforgue constructed natural compactifications of $\mathop{\rm
PGL}\nolimits_{r}^{n+1}/\mathop{\rm PGL}\nolimits_{r}$.  In fact, he
remarked that $\mathop{\rm PGL}\nolimits_{r}^{n+1}/\mathop{\rm
PGL}\nolimits_{r}$ is the quotient of $\mathop{\rm
GL}\nolimits_{r}^{n+1}/\mathop{\rm GL}\nolimits_{r}$ by the obvious free
action of the torus $\mathbb{G}_{{\rm m}}^{n+1}/\mathbb{G}_{{\rm m}}$
and that $\mathop{\rm GL}\nolimits_{r}^{n+1}/\mathop{\rm
GL}\nolimits_{r}$ may be viewed as a thin Schubert cell in the
Grassmannian variety of $r$-planes in a $r(n+1)$-dimensional vector
space.  And, more generally, he constructed natural compactifications of
all similar quotients of thin Schubert cells in the Grassmannian variety
of $r$-planes in a finite-dimensional vector space.

Let me recall that thin Schubert cells are by definition intersections
of Schubert varieties and that Israel Gelfand, Mark Goresky, Robert
MacPherson and Vera Serganova constructed natural bijections between
thin Schubert cells, matroids and certain convex polyhedra which are
called polytope matroids.

For $n=1$ and arbitrary $r$, Lafforgue's compactification of
$\mathop{\rm PGL}\nolimits_{r}^{2}/\mathop{\rm PGL} \nolimits_{r}$
coincides with the Procesi-De Concini compactification of $\mathop{\rm
PGL}\nolimits_{r}$.  It is smooth with a normal crossing divisor at
infinity.

For $n=2$ and arbitrary $r$, Lafforgue proves that his compactification
of $\mathop{\rm PGL}\nolimits_{r}^{3}/\mathop{\rm PGL} \nolimits_{r}$ is
smooth over a toric stack, and thus can be desingularized.

For $n\geq 3$ and $r\geq 3$, the geometry of Lafforgue's
compactifications is rather mysterious and not completely understood.

Gerd Faltings linked the search of good local mo\-dels for Shimura
varieties in bad characteristics to the search of smooth
compactifications of $G^{n+1}/G$ for a reductive group $G$.  He gave
another construction of Lafforgue's compactifications of $\mathop{\rm
PGL}\nolimits_{r}^{n+1}/\mathop{\rm PGL}\nolimits_{r}$ and he succeeded
in proving that Lafforgue's compactification of $\mathop{\rm
PGL}\nolimits_{r}^{n+1}/\mathop{\rm PGL}\nolimits_{r}$ is smooth for
$r=2$ and arbitrary $n$.

\section{Conclusion}

\vskip-5mm \hspace{5mm }

I hope that I gave you some idea of the depth and the technical strength
of Lafforgue's work on the Langlands correspondence for which we are now
honoring him with the Fields Medal.

\newpage

\title{\vspace*{12mm} \centerline{\Large \bf Laurent Lafforgue}\vskip 6mm}
\author{\centerline{I.H.\'{E}.S., Bures-sur-Yvette, France}\vskip 1mm}
\date{\centerline{N\'{e} le 6 novembre 1966 \`{a} Antony (Hauts-de-Seine), France} \vskip 1mm
\centerline{\large Nationalit\'{e} fran\c{c}aise}}

\maketitle

\vskip 12mm

\begin{tabular}{p{2.5cm}p{8cm}}
  1986--1990 & \'{E}l\`{e}ve \`{a} l'\'{E}cole Normale Sup\'{e}rieure de
  Paris \\
  1988--1991 & \'{E}tudiant en g\'{e}om\'{e}trie alg\'{e}brique (et en
  th\'{e}orie d'Arakelov avec Christophe Soul\'{e}) \\
  1990--1991 \newline \hspace*{4mm} puis \newline 1992--2000 & Charg\'{e} de recherche au C.N.R.S. dans
  l'\'{e}quipe ``Arithm\'{e}tique et G\'{e}om\'{e}trie Alg\'{e}brique''
  de l'Universit\'{e} Paris-Sud \\
  1991--1992 & Service militaire \`{a} l'\'{E}cole Sp\'{e}ciale
  Militaire de Saint-Cyr-Co\"{e}tquidan \\
  1993--1994 & Th\`{e}se sur les $D$-chtoucas de Drinfeld sous la
  direction de G\'{e}rard Laumon \\
  1994--2000 & Suite de l'\'{e}tude des chtoucas \\
  2000 & Professeur \`{a} l'Institut des Hautes \'{E}tudes Scientifiques
\end{tabular}

\label{lastpage}

\begin{center}
\includegraphics[scale=0.7]{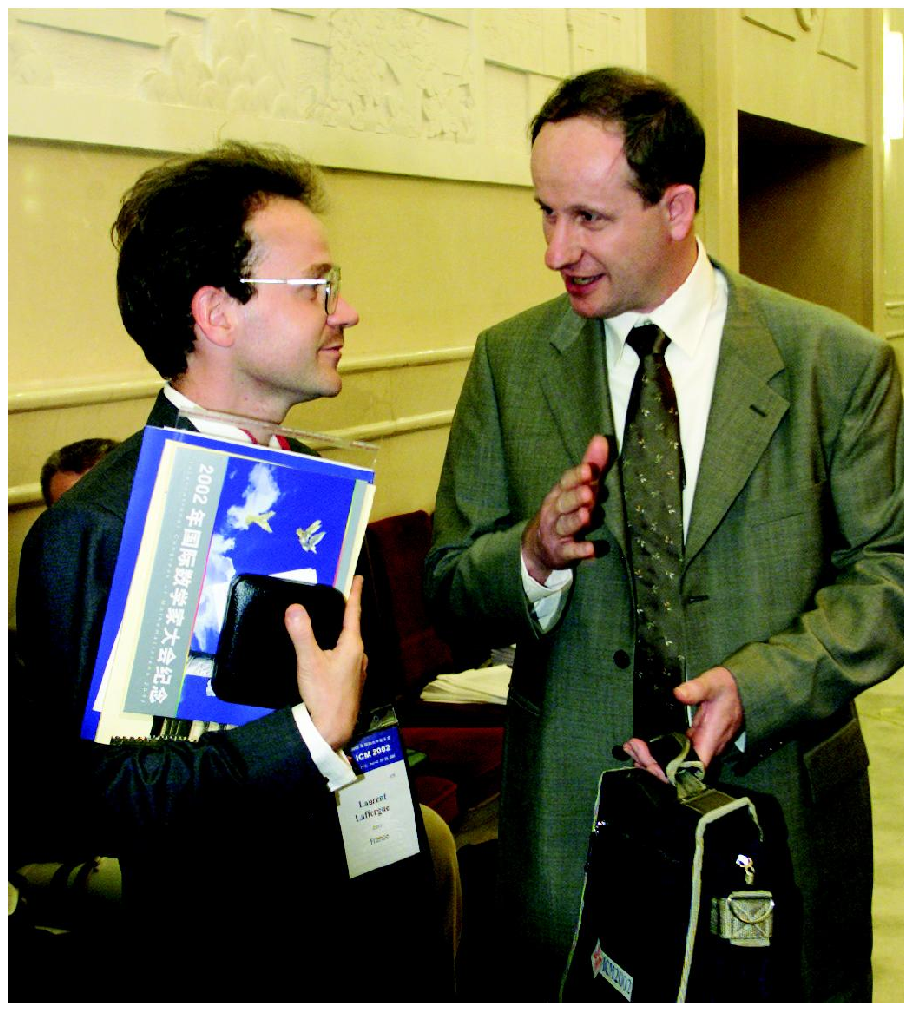}

\bigskip

L. Lafforgue (left) and G. Laumon
\end{center}

\end{document}